\documentclass[11pt]{article}

\usepackage{amsmath}
\usepackage{amssymb,amstext}
\usepackage{url}
\usepackage{array}

\newcommand{\R}{\mathbb{R}}
\newcommand{\N}{\mathbb{N}}

\newtheorem{theorem}{Theorem}[section]
\newtheorem{lemma}[theorem]{Lemma}
\newtheorem{definition}[theorem]{Definition}
\newtheorem{corollary}[theorem]{Corollary}

\newtheorem{example}[theorem]{Example}
\newtheorem{remark}[theorem]{Remark}

\DeclareMathOperator{\suppop}{supp}

\DeclareMathOperator{\rank}{rank}

\newcommand{\supp}[1]{\suppop({#1})}
\newcommand{\suchthat}{\, : \,}
\newcommand{\card}[1]{\#\,{#1}}
\newcommand{\T}{^{\rm {T}}}
\newcommand{\spark}{\mathrm {spark}}
\newcommand{\vect}{\mathrm {vec}}
\newcommand{\norm}[2]{\lVert{#1}\rVert_{#2}}

\newcommand{\abs}[1]{\lvert{#1}\rvert}

\newcommand{\define}{:=} 

\newcommand{\iprod}[2]{\langle{#1},{#2}\rangle}

\newcommand{\proof}{\par\noindent{\bf Proof}. \ignorespaces}
\newcommand{\eproof}{\space
    {\ \vbox{\hrule\hbox{\vrule height1.3ex\hskip0.8ex\vrule}\hrule}}\par}
\title{Sparse representation of solutions of
Kronecker product systems
\footnote{Supported by  \textit{Deutsche
Forschungsgemeinschaft} through the  DFG Research Center \textsc{Matheon}
\textit{Mathematics for key technologies} in Berlin.}}
\author{Sadegh Jokar
\footnote{Institut f\"ur Mathematik, TU Berlin, Str.\ des 17.\ Juni 136, 10623 Berlin, Germany.
 \texttt{$\{$jokar,mehrmann$\}$@math.tu-berlin.de}.}
\and Volker Mehrmann\footnotemark[2]
}

\date{31.07.08}

\begin{document}

\maketitle

\begin{abstract}
Three properties  of matrices: the spark, the mutual incoherence and the restricted isometry property
have recently been introduced in the context of compressed sensing.
We study these properties for matrices that are Kronecker products
and show how these properties relate to those of the factors. For the mutual incoherence we
also discuss results for sums of Kronecker products.
\end{abstract}


\begin{center}{\bf Dedicated to Shmuel Friedland  on the occasion of his 65th birthday.}
\end{center}

\section{Introduction}
In this paper we discuss the computation of sparse solutions of underdetermined linear systems
\[
Ax=b,
\]
where $A\in {\mathbb R}^{m,n} $, with $m\leq n$ is given as a Kronecker  product, i.e.
\begin{equation}\label{kp}
A= A_1\otimes A_2 \otimes \ldots \otimes A_N, \quad A_i\in \R^{m_{i},n_{i}}, i=1,\ldots,N,
\end{equation}
or as a sum of Kronecker products
\begin{equation}\label{kpsum}
A=\sum_{j=1}^{M}\,  A_1^j\otimes A_2^j \otimes \ldots \otimes A_N^j \quad A^j_i\in \R^{m_{i,j},n_{i,j}}.
\end{equation}
Since the solution is typically non-unique it is an important topic in many applications,
in particular in optimal signal recovery and in compressed sensing, see e.g.
\cite{Can06,CanR06,CanRT06,CanRT06b,CanT05,Don06a,DonE03,DonTDS06,JokP07} to find
the sparsest solution,
\begin{equation}\label{Eq:l0problem}
  \min\ \norm{x}{0},\ s.t.\ A x = b,
\end{equation}
where $\norm{x}{0}$ denotes the number of nonzero  entries of a vector $x$, see Section~\ref{sec:notation}.

In general, the  problem of finding the sparsest solution is known to be NP-hard \cite{Nat95}.
However,  in the context of compressed sensing, conditions have been derived on the size
of the \emph{support of $x$}, i.e.~the number of nonzero elements of $x$,
that allow one to compute the sparsest solution using $\ell_1$-mini\-mi\-zation via the
so called \emph{basis pursuit algorithm}~\cite{CanR06,CanRT06b,CheDS99,CohDD06a,DonE03,DonET06,DonH01}, i.e,
by computing
\begin{equation}\label{Eq:l1problem}
  \min\ \norm{x}{1},\ s.t.\ A x = b,
\end{equation}
where $\norm{x}{1}=\sum_{i} \abs{x_i}$.

Sufficient conditions for this approach to work are that some properties of the matrix  $A$
 called \emph{spark}~\cite{DonE03,Tro04b}, \emph{mutual
incoherence}~\cite{CheDS99,DonH01} or the \emph{restricted isometry property (RIP)}~\cite{Can08,CanR06,CanRT06} are restricted.
We will introduce
these properties in Section~\ref{sec:notation}.

For general matrices it is possible (though expensive) to determine the mutual incoherence, while analyzing the spark or
the restricted isometry property is  difficult. If, however, the matrix $A$ has the form (\ref{kp}) then we show
in Section \ref{sec:SparsKron} that
these properties can be easily derived from the corresponding properties of the factors.
For the mutual incoherence we can also extend these results to matrices of the form (\ref{kpsum}).

\section{Notation and preliminaries}\label{sec:notation}
For $m,n \in \N$, where $\N = \{1,2,\dots\}$, we denote by $\R^{m,n}$ the set of
real $m \times n$ matrices, by $I_n$ the $n \times n$ identity matrix, and
by $\iprod{\cdot}{\cdot}$ the Euclidean inner product in $\R^n$.
%
For $1 \leq p \leq \infty$, the $\ell_p$-norm of $x \in \R^n$ is defined by
\[
  \norm{x}{p} \define \big(\sum_{j=1}^n \abs{x_j}^p\big)^{\frac{1}{p}},
\]
with the special case
\[
\norm{x}{\infty} \define \max_{j \in \{1,\dots,n\}} \abs{x_j},
\]
if $p = \infty$.
Finally, for $x \in \R^n$, we introduce the notation
\[
\norm{x}{0} \define \card{\supp{x}},
\]
where $\supp{x} \define \{j \in \{1, \dots, n\} \suchthat x_j \neq 0\}$ is
the \emph{support} of~$x$. 
Note that $\norm{\cdot}{0}$ is not a norm, since for $\alpha \neq 0$ we have $\norm{\alpha x}{0} =\norm{x}{0}$.
We use the term \emph{k-sparse} for all vectors $x$ such that $\norm{x}{0}\leq k$.
\begin{definition}\label{Def:Kron}\cite{HorJ91,LanT85}
The \emph{Kronecker product} of $A=[a_{i,j}]\in\mathbb R^{p,q}$ and
$B=[b_{i,j}]\in \mathbb R^{r,s}$ is denoted by $A\otimes B$ and is defined to be the block matrix
\[
A \otimes B:= \left[\begin{array}{ccc} a_{1,1}B & \cdots & a_{1,q}B\\
\vdots & \ddots & \vdots\\
a_{p,1}B & \cdots & a_{p,q}B
\end{array}\right] \in \mathbb{R}^{pr,qs}.
\]
Let $C=[c_1 \cdots c_r]\in \R^{q,r}$ with columns $c_i\in \R^q$, $1 \leq i\leq r$. Then,
\[
\vect(C):=\left[
\begin{array}{c}
c_1 \\
\vdots \\
c_r
\end{array}\right]\in \mathbb R^{qr}.
\]
\end{definition}
It is well known \cite{LanT85} that the matrix equation $AXB=C$, with matrices
of appropriate dimensions, is equivalent to the linear system
$$(B^{\T} \otimes A) \vect(X)=\vect(C).$$ Furthermore, using
the perfect shuffle permutation matrices $\Pi_1,\Pi_2$, we have that
$\Pi_1 (A\otimes B)\Pi_2=B\otimes A$, see \cite{LanT85}.
%

As our first special property  we introduce the spark of a matrix.
\begin{definition}\cite{DonE03,Tro04b}\label{Def:Spark}
Let  $A=[a_1,\ldots, a_n]\in \mathbb{R}^{m,n}$, $2\leq m\leq n$ have columns $a_i$ that are normalized so that
$\norm{a_i}{2}=1$, $i=1,\ldots,n$. The \emph{spark} of $A$, denoted as  $\spark(A)$ 
is defined as the cardinality of the smallest subset of linearly dependent columns of $A$.
\end{definition}

In other words, if all $r$-dimensional subsets of column vectors of $A$ are linearly independent,
but there exists a subset of $r+1$ columns that are linearly dependent, then $\spark(A)=r+1$.
For convenience, if $m=n=1$, we define $\spark(A):=1$, and in the case where $m=n\geq 2$ and $A$ is invertible, we
set $\spark(A):=n+1$.
 In general the $\spark$ and the $\rank$ of a matrix $A\in \mathbb{R}^{m,n}$
with $m\geq 2$, are related via
\[
2 \leq \spark(A)\leq \rank(A)+1.
\]
\begin{example}\label{ex1}
If
\[
A=\left[\begin{array}{rrrr}
1 & 0 & 1 & 1 \\
0 & 1 & 1 & -1
\end{array}\right],
\]
then $\spark(A)=\rank(A)+1=3$. On the other hand, if
\[
A=\left[\begin{array}{rrrr}
1 & 0 & 1 & -1 \\
0 & 1 & 1 & 0
\end{array}\right],
\]
then $\spark(A)=2$.
\end{example}
The quantity $\spark(A)$ can be used to derive necessary and
sufficient conditions for the existence of sparse solutions.
\begin{lemma}\label{Lem:SparseSpark}\cite{DonE03,GriN03}
Consider the linear system $Ax=b$ with $A \in \mathbb R^{m,n}$, $m\leq n$.
A necessary and sufficient condition for the linear system $A x=b$ to have a
unique $k$-sparse solution $x$ is that $k \leq \spark(A)/2$.
\end{lemma}

The second property that we study is the mutual incoherence.
\begin{definition}\label{Def:Mut1}\cite{DonH01}
Let $A =[a_1,\ldots, a_n]\in \mathbb R^{m,n}$, $m\leq n$ have columns $a_i$
that are normalized so that $\norm{a_i}{2}=1$, $i=1,\ldots,n$. Then the
{\emph mutual incoherence ${\mathcal M}(A)$} is defined by
\[
{\mathcal M}(A):= \max_{i\neq j} |\langle a_i,a_j\rangle|=\max_{i\neq j} |(A^{\T} A)_{i,j}|.
\]
\end{definition}
Note that, since the columns of $A$ are normalized,
by the triangle inequality we always have ${\mathcal M}(A)\leq 1$.
On the other hand, if $A$ has orthonormal columns, then ${\mathcal M}(A)=0$.

We have the following lower bound for ${\mathcal M}(A)$.
\begin{lemma}\label{lem:UpBoundMI}\cite{StrH03}
Suppose that $A\in \mathbb{R}^{m,n}$, $m\leq n $ has  columns $a_i$
that are normalized so that $\norm{a_i}{2}=1$, $i=1,\ldots,n$ and suppose further that $A$ has full row rank. Then
\[
{\mathcal M} (A)\geq \sqrt{\frac{n-m}{m(n-1)}}.
\]
\end{lemma}
%
The following lemma relates the sparsest solution as defined in (\ref{Eq:l0problem}) and the
$\ell_1$-solution as defined in (\ref{Eq:l1problem}) of the linear equation $Ax=b$ in terms of the mutual incoherence
of a matrix $A$.
\begin{lemma}\cite{DonE03,Fuc05}\label{Lem:MIl0}
Suppose that $A\in \mathbb{R}^{m,n}$, $m\leq n $ has  columns $a_i$
that are normalized so that $\norm{a_i}{2}=1$, $i=1,\ldots,n$.
 If $b$ is a vector such that the equation $Ax=b$ has a solution satisfying
\[
\norm{x}{0}< \frac{1+ \frac{1}{{\mathcal M}(A)} } {2},
\]
 then the $\ell_1$-norm minimal solution in (\ref{Eq:l1problem})
coincides with the $\ell_0$-minimal solution in (\ref{Eq:l0problem}).
\end{lemma}
\begin{remark}
Consider matrices of the form  $A =[\Phi\; \Psi]$, where $\Phi$ and $\Psi$ have orthonormal columns.
If the sparsest solution $x$ of $Ax=b$ satisfies
$$\norm{x}{0}< \frac{\sqrt{2}-\frac{1}{2}}{{\mathcal M}(A)},$$
then it has been shown in \cite{ElaB02} that the solutions of the
$\ell_1$-norm mini\-mi\-zation problem and $\ell_0$-norm mini\-mi\-zation problem coincide.
\end{remark}

The third quantity that is important in the context of sparse recovery and
compressed sensing is the \emph{restricted isometry property}.
\begin{definition}\label{Def:RIP}~\cite{Can08,CanR06,CanRT06,CanRT06b}
Let $A =[a_1,\ldots, a_n]\in \mathbb R^{m,n}$, $m\leq n$ have columns $a_i$
that are normalized so that $\norm{a_i}{2}=1$, $i=1,\ldots,n$.

The \emph{$k$-restricted isometry constant of $A$} is the smallest number $\delta_k$ such that
\[
(1-\delta_k)\norm{x}{2}^2 \leq \norm{A x}{2}^2 \leq (1+\delta_k) \norm{x}{2}^2
\]
for all $x\in \mathbb{R}^n$ with $\norm{x}{0} \leq k$.
\end{definition}
The $k$-restricted isometry property
requires that every set of columns of cardinality less than or equal to $k$ approximately (with an error $\delta_k$)
behaves like an
orthonormal basis. 

The following lemma gives the relation between the sparsest solution (as defined in (\ref{Eq:l0problem})) of a
linear system  $Ax=b$  and the
$\ell_1$-solution as defined in (\ref{Eq:l1problem}) in terms of the $k$-restricted isometry constant.
\begin{lemma}\cite{Can08}\label{Lem:RIPSp}
Let $A =[a_1,\ldots, a_n]\in \mathbb R^{m,n}$, $m\leq n$ have columns $a_i$
that are normalized so that $\norm{a_i}{2}=1$, $i=1,\ldots,n$.

Suppose that
\[
\delta_{2k} < \sqrt{2}-1.
\]
Then for all $k$-sparse solution vectors $x$ of $Ax=b$
 the solution of (\ref{Eq:l1problem}) is equal to the solution of (\ref{Eq:l0problem}).
\end{lemma}

After introducing the concepts of spark, mutual incoherence and $k$-restricted isometry property,
in the next section we analyze these concepts for Kronecker product matrices.

\section{Sparse representation and Kronecker Products of Matrices}\label{sec:SparsKron}
In this section we study sparse solutions for linear system $Ax=b$, where the matrix $A$ is given as a
Kronecker product (\ref{kp}).

Our first result characterizes $\spark(A \otimes B)$ in terms of $\spark(A)$ and $\spark(B)$.
Note  that if $A,B$ have normalized columns then $A\otimes B$ has normalized columns as well.
\begin{theorem}\label{Thm:SparkKron}
Let $A=[a_1,\ldots,a_q]\in \R^{p,q}$ and $B=[b_1,\ldots,b_s] \in \R^{r,s}$
be rank-deficient matrices with normalized columns,
i.e., $\norm{a_i}{2}=1$, $i=1,\ldots,q$, $\norm{b_i}{2}=1$, $i=1,\ldots,s$. Then
\begin{equation}\label{Eq:SparkKron}
\spark(A \otimes B)=\spark(B \otimes A)= \min\{\spark(A),\spark(B)\}.
\end{equation}
If $A$ is an invertible matrix and $B$ is rank-deficient matrix, then
\begin{equation}\label{Eq:SparkKronAInv}
\spark(A \otimes B)=\spark(B).
\end{equation}
If both $A$ and $B$ are square and invertible then
\[
\spark(A \otimes B)=(\spark(A)-1)(\spark(B)-1)+1=qs+1.
\]
\end{theorem}
\proof
Using the fact that $(B\otimes A) \vect(X)=\Pi_1 (A \otimes B) \Pi_2 \vect(X)$ and
 $\norm{\vect(X)}{0}=\norm{\Pi_2 \vect(X)}{0}$, we have $\spark(A \otimes B)=\spark(B \otimes A)$.

Consider first the case that $A$ and $B$ are rank-deficient.
By the definition of $\spark(B)$, there exists a vector
$y \in \R^{s}$ with $\norm{y}{0}=\spark(B)$ such that $B y=0$. With
$$
\hat{X}=[\begin{array}{c c c c} y & 0 & \cdots & 0 \end{array}],
$$
we have that $(A \otimes B) \vect(\hat{X})=0$ and $\norm{\vect(\hat{X})}{0}=\norm{y}{0}=\spark(B)$. This means that
$\spark(A\otimes B) \leq \spark(B)$. Using that $\spark(A \otimes B)=\spark(B \otimes A)$ and that
also $A$ is rank-deficient, we can apply the same argument as before and get
 $\spark(A \otimes B) \leq \spark(A)$. Therefore,
\begin{equation}\label{Eq:MinSpark}
\spark(A \otimes B)\leq \min\{\spark(A),\spark(B)\}.
\end{equation}

Let $C=A\otimes B$, then every column of $C$ has the form $c_j=a_{u_j}\otimes b_{v_j}$.
To prove equality in (\ref{Eq:SparkKron}), we assume w.l.o.g. that
\begin{equation}\label{Eq:AssumSparkAB}
\spark(B)\leq \spark(A).
\end{equation}
Then by (\ref{Eq:MinSpark}) we have
$
\spark(A \otimes B)\leq \spark(B).
$
Suppose now that
\begin{equation}\label{Eq:ContrSpark}
\spark(A \otimes B)=\ell<\spark(B).
\end{equation}
This implies, in particular, that any set of $\ell$ columns of $B$ is linearly independent, while
there exist scalars $\lambda_1,\ldots,\lambda_\ell$ not all $0$ and indices
$u_1,\ldots ,u_\ell$ where $u_i\neq u_j$ for all $i\neq j$, and  $v_1,\ldots,v_\ell$ such that
\[
\sum_{j=1}^\ell (a_{u_j}\otimes b_{v_j})\lambda_j=\sum_{j=1}^\ell (\lambda_j a_{u_j})\otimes b_{v_j}=0.
\]
In this sum there may occur repeated copies of vectors $b_j$, so without loss of generality we may assume the indices
$v_i$ are numbered so that
\[
\underbrace{v_1=\cdots=v_{k_1}}_{g_1}< \underbrace{v_{k_1+1}=\cdots=v_{k_2}}_{g_2}
< \cdots < \underbrace{v_{k_{t-1}+1}=\cdots=v_{k_{t}}}_{g_{t}}.
\]
Therefore, we have
\begin{equation}\label{Eq:SparkLinear1}
(\sum_{j=1}^{k_1}\lambda_j a_{u_j})\otimes b_{g_1} + (\sum_{j=k_1+1}^{k_2}
\lambda_j a_{u_j})\otimes b_{g_2}+\cdots+(\sum_{j=k_{t-1}+1}^{k_{t}}  \lambda_j a_{u_j}) b_{g_{t}}=0,
\end{equation}
where $k_{t}=\ell$. Since $b_{g_1},\ldots,b_{g_{t}}$ are linearly independent,
it follows that for all $1 \leq i \leq t$ we have
\[\sum_{j=k_{i-1}+1}^{k_i}\lambda_j a_{u_j}=0,\]
where $k_0=0$. This contradicts the assumption in (\ref{Eq:AssumSparkAB}) that
\[
\ell<\spark(B) \leq \spark(A),
\]
because the $u_j$ are pairwise distinct and at least one of the coefficients $\lambda_j$ is nonzero.

Now suppose that $A$ is invertible and $B$ is rank-deficient. Then with the same argument as above, we have
$\spark(A\otimes B)\leq \spark(B)$. Let $X=[x_1,\ldots, x_q]\neq 0$, such that
$\norm{\vect(X)}{0}=\spark(A\otimes B)$ and $(A\otimes B) \vect(X)=0$. This implies that $BXA^{\T}=0$,
 and, since $A$ is invertible we have $BX=0$, while on the other hand $X\neq 0$.
Thus there exists at least one index $i$ such that $x_i\neq 0$ and $Bx_i=0$. Hence,
\[
\spark(B)\leq \norm{x_i}{0} \leq  \norm{\vect(X)}{0}=\spark(A\otimes B),
\]
and therefore $\spark(A\otimes B)=\spark(B)$.

For the case where both $A$ and $B$ are invertible, $A\otimes B$ is invertible as well, see \cite{LanT85}.
 Therefore,
\[\spark(A\otimes B)=\rank(A\otimes B)+1=qs+1.
\]
\eproof

We immediately have the following corollary of Theorem \ref{Thm:SparkKron}.
\begin{corollary}\label{Cor:SparkKron}
Consider rank-deficient matrices $\{A_i\}_{i=1}^N$ with normalized columns. Then
\[
\spark(A_1\otimes \ldots \otimes A_N)=\displaystyle\min_{1 \leq i\leq N} \{\spark(A_i)\}.
\]
\end{corollary}
By combining Lemma \ref{Lem:SparseSpark} and Corollary \ref{Cor:SparkKron} we get the following Corollary:
\begin{corollary}\label{Cor:SparkKronl0}
Consider a linear system  $(A_1\otimes \ldots \otimes A_N)x=b$ with rank-deficient matrices
$A_i\in \R^{p_i,q_i}$ that have normalized columns. A necessary and sufficient condition for this linear system
to have a unique $k$-sparse solution $x$ is that
\[k \leq \frac{\displaystyle\min_{1 \leq i \leq N}\{\spark(A_i)\}}{2}.\]
\end{corollary}
\begin{remark}
Corollary \ref{Cor:SparkKronl0} is saying that if one of the matrices $A_j$ has small
$\spark$ then we can only uniquely recover vectors of the sparsity up to
$\spark(A_j)/2$ in the linear system $(A_1\otimes\ldots\otimes A_N)x=b$.
\end{remark}

Similar to the analysis of $\spark(A\otimes B)$, we can also obtain
an estimate of ${\mathcal M}(A \otimes B)$ in terms of ${\mathcal M}(A)$ and
${\mathcal M}(B)$.
\begin{theorem}\label{Thm:MutKron}
Consider matrices $A=[a_1,\ldots,a_{n_1}]\in \R^{m_1,n_1}$ and
 $B=[b_1,\ldots,b_{n_2}]\in \R^{m_2,n_2}$ with normalized columns.
Then
\[
{\mathcal M}(A \otimes B)=\max\{{\mathcal M}(A),{\mathcal M}(B)\}.
\]
\end{theorem}
\proof
 Suppose that $C=A\otimes B$ and $C=[c_1 \cdots c_n]$, where $c_i\in \mathbb{R}^m$, $m=m_1 m_2$ and $n=n_1 n_2$.
 Then we have ${\mathcal M}(C)=\max_{i\neq j} \abs{\langle c_i, c_j\rangle}$. Since
 $c_i=a_p\otimes b_q$ and $c_j=a_r \otimes b_s$ for some $p,q,r,s$,
using properties of the Kronecker product~\cite{LanT85}, we have
\begin{equation}\label{Eq:Kron1}
\langle c_i, c_j \rangle=\langle a_p\otimes b_q, a_r \otimes b_s
\rangle=\langle a_p, a_r\rangle \cdot \langle b_q, b_s\rangle.
\end{equation}
By Definition~\ref{Def:Mut1} and (\ref{Eq:Kron1}) we then have
\begin{eqnarray}\label{Eq:Kron2}
\begin{array}{l l}
{\mathcal M}(C)={\mathcal M}(A\otimes B) & =\displaystyle\max_{\substack{ p,q,r,s \\ (p,q) \neq (r,s)}} \abs{\langle
a_p, a_r\rangle \cdot \langle b_q, b_s\rangle}\\
& =\displaystyle\max_{\substack{ p,q,r,s \\ p\neq r,q\neq s}}
\{\abs{\langle a_p, a_r\rangle \cdot \langle b_q, b_s\rangle} , \abs{\langle a_p,
a_r\rangle} , \abs{\langle b_q, b_s\rangle} \}.
\end{array}
\end{eqnarray}
On the other hand, since the matrices $A$ and $B$ have normalized columns, we have
\[\abs{\langle a_p, a_r\rangle \cdot \langle b_q, b_s\rangle} \leq \abs{\langle a_p, a_r\rangle},\]
and similarly
\[\abs{\langle a_p, a_r\rangle \cdot \langle b_q, b_s\rangle} \leq \abs{\langle b_q, b_s\rangle}.\]
Therefore, from (\ref{Eq:Kron2}) we have
\begin{eqnarray*}
\begin{array}{l l}
{\mathcal M}(A\otimes B) & =\displaystyle\max_{\substack{p,q,r,s \\ p\neq r,q\neq s}}
\{ \abs{\langle a_p, a_r\rangle}, \abs{\langle b_q, b_s\rangle} \} \\
& =\max\{\displaystyle\max_{p\neq r}\abs{\langle a_p, a_r\rangle},
 \displaystyle\max_{q\neq s}\abs{\langle b_q, b_s\rangle}\} \\
& =\max\{{\mathcal M}(A),{\mathcal M}(B)\}.
\end{array}
\end{eqnarray*}
\eproof
A direct consequence of this theorem is the following Corollary.
\begin{corollary}\label{Cor:Kron}
Consider matrices $\{A_i\}_{i=1}^N$ with normalized columns and let $A=A_1\otimes\ldots \otimes A_N$.
Then,
\[
{\mathcal M}(A)=\max_{1\leq i \leq n}{\mathcal M}(A_i).
\]
\end{corollary}
 Corollary \ref{Cor:Kron} shows that if one of the matrices $A_i$ has a large mutual incoherence,
 then it will dominate the mutual incoherence of $A$, regardless of all the other factors in
 the Kronecker product.

We also have a result that relates the $k$-restricted isometry constant
of $\delta_k^{A \otimes B}$ to those of $\delta_k^{A}$ and $\delta_k^{B}$.
\begin{theorem}\label{Thm:RIPKron}
Let
$A\in \R^{p,q}$ and $B\in \R^{r,s}$ have normalized columns. Then
\begin{equation}\label{eq:RIPKron}
\delta^{A \otimes B}_k=\delta^{B \otimes A}_k \geq \max\{\delta^{A}_k,\delta^{B}_k\}.
\end{equation}
\end{theorem}
\proof Using the fact that $B \otimes A=\Pi_1 (A \otimes B) \Pi_2$, where $\Pi_1$ and $\Pi_2$ are permutation matrices we have
\[
\norm{\vect(X)}{2}^2=\norm{\Pi_2 \vect(X)}{2}^2,
\]
and
\[
\norm{(B \otimes A) \vect(X)}{2}^2=
\norm{\Pi_1 (A \otimes B) \Pi_2 \vect(X)}{2}^2 = \norm{(A \otimes B) (\Pi_2 \vect(X))}{2}^2.
\]
Therefore $\delta^{A \otimes B}_k=\delta^{B \otimes A}_k$. To prove the assertion,
it is sufficient to prove that $\delta^{A\otimes B}_k\geq \delta^{B}_k$,
 the proof that $\delta^{A\otimes B}_k\geq \delta^{A}_k$
follows analogously. We know that $\delta_k^{B}$ is the smallest constant such that, for all $x$ with
$\norm{x}{0}\leq k$, we have
\[
(1-\delta_k^{B})\norm{x}{2}^2 \leq \norm{B x}{2}^2 \leq (1+\delta_k^{B})\norm{x}{2}^2.
\]
For any $x$ with $\norm{x}{0}\leq k$, we can construct the
 matrix $X=[\begin{array}{l l l l } x & 0 & \cdots & 0 \end{array}]$, with
$\norm{\vect(X)}{0} \leq k$. Since $A$ has normalized columns, we have
\begin{equation}\label{Eq:Kroneq}
\norm{(A \otimes B) (\vect(X))}{2}^2=\sum_{i=1}^p a_{i,1}^2 \norm{B x}{2}^2=\norm{B x}{2}^2,
\end{equation}
and
\begin{equation}\label{Eq:veceq}
\norm{\vect(X)}{2}^2=\norm{x}{2}^2.
\end{equation}
On the other hand $\delta_k^{A \otimes B}$ is the smallest constant such that
\[
(1-\delta_k^{A \otimes B})\norm{\vect(X)}{2}^2 \leq \norm{(A \otimes B) (\vect(X))}{2}^2
\leq (1+\delta_k^{A \otimes B})\norm{\vect(X)}{2}^2,
\]
and for the special class of $k$-sparse vectors $\vect(X)$
from (\ref{Eq:Kroneq}) and (\ref{Eq:veceq}) we have
\[
(1-\delta_k^{B})\norm{\vect(X)}{2}^2 \leq \norm{(A \otimes B) (\vect(X))}{2}^2
 \leq (1+\delta_k^{B})\norm{\vect(X)}{2}^2,
\]
where $\delta_k^{B}$ is the smallest constant for this special class of $k$-sparse vectors.
Therefore, for general $k$-sparse vectors, we have
$$\delta_k^{A \otimes B} \geq \delta_k^{B}.$$
\eproof

\begin{remark}
 Note that for $k=2$, equality holds in (\ref{eq:RIPKron}), since 
for a given normalized matrix $A$, we have $\delta_2^{A}={\mathcal M}(A)$.  Therefore, by Theorem~\ref{Thm:MutKron}
it follows that
 \[
\delta^{A \otimes B}_2=\max\{\delta^{A}_2,\delta^{B}_2\}.
\]
 For $k\geq 3$, however the inequality may be strict.
For example if $A=[I_2 \ H_2]$ and $B=[I_4 \ H_4]$, where $I_n$ is the identity matrix of order $n$
and $H_n$ is the \emph{Hadamard} matrix of order $n$, see e.g.~\cite{Had1893,Hor07},
then $\delta_3^{A \otimes B}=1.0545 > \max\{1,\frac{1}{\sqrt{2}}\}$.
Here the $k$-restricted isometry constants of these matrices were calculated using the singular value
decomposition for all submatrices consisting of $3$ columns.
\end{remark}

We have the obvious corollary.
\begin{corollary}\label{Cor:RIPKron}
Suppose that matrices $A_i$ for $i=1,\ldots,N$ have normalized columns. Then
\[
\delta^{A_1 \otimes \ldots \otimes A_N}_k \geq \max_{1 \leq i\leq N} \{\delta^{A_i}_k\}.
\]
\end{corollary}

According to Lemma \ref{Lem:RIPSp}, if the restricted isometry constant $\delta_{2k}$
is small enough ($\delta_{2k}<\sqrt{2}-1$), then one can recover all $k$-sparse solutions using $\ell_1$-mini\-mi\-zation.
On the other hand, Corollary \ref{Cor:RIPKron} implies that if the $k$-restricted isometry constant
$\delta_k$ of $A$ is small (for example less than $1/2$), then $A$ can not be written as a Kronecker product of matrices
$A_i$ with smaller sizes.

\section{Sums of Kronecker products}\label{sec:sumkp}
In many applications, in particular in finite difference or finite element discretizations of partial differential equations
in more than one space dimension~\cite{Str07}, linear systems with matrices that are sums of Kronecker products arise.

It is then an obvious question whether the spark, the mutual incoherence and the $k$-restricted isometry
property for sums of Kronecker products can be related to that of the summands.

Unfortunately, in general we do not have a nice relation between $\spark(A+B)$ and  $\spark(A)$, $\spark(B)$.

\begin{example}\label{exsp}
Let $E_n$ denote the $n\times n$ matrix of all ones.
If
\[ A=\left[\begin{array}{ll} I_{5} & E_{5} \end{array}\right]
\otimes \left[\begin{array}{rr} 1 & 1 \\ 1 &  -1\end{array}\right],
\]
and
\[ B=\left[\begin{array}{ll} [1 \ 2 \ 3 \ 4 \ 5]^{\T} [1 \ 2 \ 3 \ 4 \ 5]
 & I_{5} \end{array}\right] \otimes \left[\begin{array}{rr} 1 & 1 \\ 0 &  -1\end{array}\right],
\]
then $5=\spark(A+B) > \spark(A)+\spark(B) =2+2$.

  On the other hand if $A=I_{2}\otimes I_{2}$ and
$A+B=\frac{1}{2} (E_{2}\otimes E_{2})$ then $2=\spark(A+B)<\spark(A)+\spark(B)=5+5$.
\end{example}

For the mutual incoherence the situation is better. We introduce the following concept of
diagonal and off-diagonal mutual incoherence.
\begin{definition}\label{Def:Mut2}
Suppose that $A=[a_1,\ldots, a_n],\ B=[b_1,\ldots, b_n]\in \mathbb{R}^{m,n}$, $m\leq n $,
have normalized columns. Then the \emph{off-diagonal mutual
incoherence ${\mathcal M}_{OD}(A,B)$ of $A$ and $B$} is defined via
\[
{\mathcal M}_{OD}(A,B):= \max_{i\neq j} |\langle a_i, b_j\rangle|
\]
and the \emph{diagonal mutual incoherence ${\mathcal M}_{D}(A,B)$  of $A$ and $B$} is defined via
\[
{\mathcal M}_{D}(A,B):= \max_{i} |\langle a_i,b_i \rangle|.
\]
\end{definition}

\begin{remark}
Note that in Definition~\ref{Def:Mut2} the order of the columns is
important. Note further that in the special case that $A=B$ we have
${\mathcal M}_{OD}(A,A)={\mathcal M}(A)$ and ${\mathcal M}_{D}(A,A)=1$.
\end{remark}
%

Then we have the following theorem.
\begin{theorem}\label{Summutual1}
Let $A=[a_1,\ldots, a_n], B=[b_1,\ldots, b_n]\in \mathbb R^{m,n}$
be matrices with normalized columns and suppose that ${\mathcal M}_{D}(A,B)\neq 1$.  Then,
%
\begin{equation}\label{Mutualsum1}
{\mathcal M}(A+B)\leq \frac{{\mathcal M}(A)+2{\mathcal M}_{OD}(A,B)+{\mathcal M}(B)}{2(1-{\mathcal M}_{D}(A,B))}.
\end{equation}
\end{theorem}
\proof
For $i\neq j$, by the triangle inequality we have that
\begin{eqnarray}\label{Triag1}\begin{array}{ll}
\abs{\langle a_i+b_i,a_j+b_j \rangle} & \leq \abs{\langle a_i,a_j \rangle}+
\abs{\langle a_i,b_j \rangle}+\abs{\langle b_i,a_j \rangle}+\abs{\langle b_i,b_j \rangle} \\ \\
& \leq {\mathcal M}(A)+2{\mathcal M}_{OD}(A,B)+{\mathcal M}(B)
\end{array}
\end{eqnarray}
and
\begin{equation}\label{Triag2}
\norm{a_i+b_i}{2}^2=2+2 \langle a_i,b_i \rangle \geq 2(1-\abs{ \langle a_i,b_i \rangle}) \geq 2(1-{\mathcal M}_D(A,B)).
\end{equation}
Combining (\ref{Triag1}) and (\ref{Triag2}), we get
\[
{\mathcal M}(A+B)=\max_{i\neq j} \frac{\abs{\langle a_i+b_i,a_j+b_j \rangle}}{\norm{a_i+b_i}{2}\norm{a_j+b_j}{2}}
\leq \frac{{\mathcal M}(A)+2{\mathcal M}_{OD}(A,B)+{\mathcal M}(B)}{2(1-{\mathcal M}_{D}(A,B))}.
\]
\eproof
Note that the inequality (\ref{Mutualsum1}) also holds if  ${\mathcal M}_{D}(A,B)=1$, if we define the right side
to be infinite in this case.
\begin{remark}\label{remsharp}
The bound in Theorem~\ref{Summutual1} is sharp. For example if
\begin{eqnarray*}
A=\left[\begin{array}{l l}
0 & 1 \\
1 & 0
\end{array}\right], \quad B=I_2
\end{eqnarray*}
then
\[{\mathcal M}(A+B)=1, \quad {\mathcal M}(A)={\mathcal M}(B)={\mathcal M}_D(A,B)=0\]
and
\[{\mathcal M}_{OD}(A,B)=1.
\]
\end{remark}
Theorem~\ref{Summutual1} immediately extends to more than one summand.
\begin{corollary}\label{Summutual2}
Consider matrices  $A^i\in \mathbb R^{m,n}, \ 1 \leq i \leq M$ with normalized columns. If
\[M-2\sum_{1 \leq i<j\leq M} {\mathcal M}_{D}(A^i,A^j) > 0,\]
then
\[
{\mathcal M}(\sum_{i=1}^M A^i)\leq \frac{\sum_{i=1}^M {\mathcal M}(A^i)+2
\sum_{1 \leq i<j\leq M} {\mathcal M}_{OD}(A^i,A^j)}{M-2\sum_{1 \leq i<j\leq
M} {\mathcal M}_{D}(A^i,A^j)}.
\]
\end{corollary}
%

In order to apply these results to sums of Kronecker products of
the form $A=\sum_{j=1}^M A^j_1\otimes \ldots \otimes A^j_N$,
we introduce the abbreviation
\[ {\mathcal U}(\sum_{j=1}^M A^j ):=\frac{\sum_{j=1}^M {\mathcal M}(A^j)+2 \sum_{1 \leq i<j\leq M}
 {\mathcal M}_{OD}(A^i,A^j)}{M-2\sum_{1 \leq i<j\leq M} {\mathcal M}_{D}(A^i,A^j)}.
\]
We have the following Corollary.
\begin{corollary}\label{Thm:MutKronl0} Consider a linear system of the form
\[
(\sum_{j=1}^M A^j_1\otimes \ldots\otimes A^j_N)x=b,
\]
where the matrices $A_i^j$ are of appropriate dimensions and have normalized columns.
Suppose that there exists a solution $x$ with the sparsity
\[\norm{x}{0} < \frac{1}{2} \left(1+\frac{1}{{\mathcal U}(\sum_{j=1}^M A^j_1\otimes \ldots\otimes A^j_N}\right).
\]
Then this is the  unique solution with this sparsity which can be recovered using $\ell_1$-mini\-mi\-zation as defined in (\ref{Eq:l1problem}).
\end{corollary}
\proof By applying Lemma \ref{Lem:MIl0} and Corollary \ref{Summutual2},  we have that
\[\norm{x}{0} < \frac{1}{2} \left(1+\frac{1}{{\mathcal U}(\sum_{j=1}^M A^j_1\otimes \ldots\otimes A^j_N}\right)
\]
implies that
\[\norm{x}{0} < \frac{1}{2} \left(1+\frac{1}{{\mathcal M}(\sum_{j=1}^M A^j_1\otimes \ldots\otimes A^j_N}\right)
\]
and therefore by Lemma \ref{Lem:MIl0} the sparse solution is unique.
\eproof

\begin{example}
Consider a sum of Kronecker products
$C=I\otimes A+A \otimes I$ as they for example arise in the finite difference approximation
of boundary value problems for $2D$ elliptic PDEs. Then it is easy to see that
\[
{\mathcal M}(C) \leq {\mathcal U}(C)=\frac{2 {\mathcal M}(A)+ 2
{\mathcal M}_{OD}(I\otimes A, A\otimes I)}{2-2{\mathcal M}_{OD}(I\otimes A, A\otimes I)}.
\]
Especially,  if  $A$ is a $1D$ finite difference matrix, e.g.
\[
A=\left[\begin{array}{rrr} 2 & -1 & 0 \\ -1 & 2 & -1 \\ 0 & -1 & 2 \end{array}\right],\]
then we have
\[
0.4237={\mathcal M}(C) \leq {\mathcal U}(C)=0.5615.
\]
\end{example}

For the $k$-restricted isometry property it is an open problem to establish relationships
between that of a sum of Kronecker products and the summands.
%
\section{Conclusion}
We have analyzed the recently introduced concepts of the spark, the mutual incoherence and the $k$-restricted
isometry property of matrix in Kronecker product form to that of the Kronecker factors.

\section*{Acknowledgment}
We thank J. Gagelman, O. Holtz, M. Pfetsch, C. Van Loan,  Z. Xu and H. Yserentant for
 fruitful discussions on this topic.
\bibliographystyle{mod_siam}
\bibliography{JokM08}

\begin{thebibliography}{10}

\bibitem{Can06}
\textsc{E.~J. Cand\'{e}s}, \emph{Compressive sampling}, in Proc. International
  Congress of Mathematics, Madrid, Spain, 2006, pp.~1433--1452.

\bibitem{Can08}
\textsc{E.~J. Cand\`{e}s}, \emph{The restricted isometry property and its
  implications for compressed sensing}.
\newblock Technical Report, California Institute of Technology, 2008.

\bibitem{CanR06}
\textsc{E.~J. Cand\`{e}s and J.~Romberg}, \emph{Quantitative robust uncertainty
  principles and optimally sparse decompositions}, Found. Comput. Math.
  \textbf{6}, no.~2 (2006), pp.~227--254.

\bibitem{CanRT06}
\textsc{E.~J. Cand\`{e}s, J.~Romberg, and T.~Tao}, \emph{Robust uncertainty
  principles: Exact signal reconstruction from highly incomplete frequency
  information}, IEEE Trans. Inform. Theory \textbf{52}, no.~2 (2006),
  pp.~489--509.

\bibitem{CanRT06b}
\textsc{E.~J. Cand\`{e}s, J.~Romberg, and T.~Tao}, \emph{Stable signal recovery
  from incomplete and inaccurate measurements}, Comm. Pure Appl. Math.
  \textbf{59}, no.~8 (2006), pp.~1207--1223.

\bibitem{CanT05}
\textsc{E.~J. Cand\`{e}s and T.~Tao}, \emph{Decoding by linear programming},
  IEEE Trans. Inform. Theory \textbf{51}, no.~12 (2005), pp.~4203--4215.

\bibitem{CheDS99}
\textsc{S.~S. Chen, D.~L. Donoho, and M.~A. Saunders}, \emph{Atomic
  decmposition by basis pursuit}, SIAM J. Sci. Comput. \textbf{20}, no.~1
  (1999), pp.~33--61.

\bibitem{CohDD06a}
\textsc{A.~Cohen, W.~Dahmen, and R.~DeVore}, \emph{Compressed sensing and best
  k-term approximation}.
\newblock Preprint, 2006.

\bibitem{Don06a}
\textsc{D.~L. Donoho}, \emph{Compressed sensing}, IEEE Trans. Inform. Theory
  \textbf{52}, no.~4 (2006), pp.~1289--1306.

\bibitem{DonE03}
\textsc{D.~L. Donoho and M.~Elad}, \emph{Optimally sparse representation in
  general (non\-orthogonal) dictionaries via $\ell^1$ minimization}, Proc.
  Natl. Acad. Sci. USA \textbf{100}, no.~5 (2003), pp.~2197--2202.

\bibitem{DonET06}
\textsc{D.~L. Donoho, M.~Elad, and V.~Temlyakov}, \emph{Stable recovery of
  sparse overcomplete representations in the presence of noise}, IEEE Trans.
  Inform. Theory \textbf{52}, no.~1 (2006), pp.~6--18.

\bibitem{DonH01}
\textsc{D.~L. Donoho and X.~Huo}, \emph{Uncertainty principles and ideal atomic
  decomposition}, IEEE Trans. Inform. Theory \textbf{47}, no.~7 (2001),
  pp.~2845--2862.

\bibitem{DonTDS06}
\textsc{D.~L. Donoho, Y.~Tsaig, I.~Drori, and J.-L. Starck}, \emph{Sparse
  solution of underdetermined linear equations by stagewise orthogonal matching
  pursuit}, Tech. Report 2006-02, Stanford, Department of Statistics, 2006.

\bibitem{ElaB02}
\textsc{M.~Elad and A.~M. Bruckstein}, \emph{A generalized uncertainty
  principle and sparse representation in pairs of bases}, IEEE Trans. Inform.
  Theory \textbf{48}, no.~9 (2002), pp.~2558--2567.

\bibitem{Fuc05}
\textsc{J.~J. Fuchs}, \emph{Recovery of exact sparse representations in the
  presence of bounded noise}, IEEE Trans. Inform. Theory \textbf{51}, no.~10
  (2005), pp.~1601--1608.

\bibitem{GriN03}
\textsc{R.~Gribonval and M.~Nielsen}, \emph{Sparse representations in unions of
  bases}, IEEE Trans. Inform. Theory \textbf{49}, no.~12 (2003),
  pp.~3320--3325.

\bibitem{Had1893}
\textsc{J.~Hadamard}, \emph{R\'{e}solution d'une question relative aux
  d\'{e}terminants}, Bull. Des Sciences Math. \textbf{17}, no.~2 (1893),
  pp.~240--246.

\bibitem{Hor07}
\textsc{K.~J. Horadam}, \emph{Hadamard Matrices and Their Applications},
  Princton University Press, Princeton, NJ, 2007.

\bibitem{HorJ91}
\textsc{R.~Horn and C.~R. Johnson}, \emph{Topics in Matrix Analysis}, Cambridge
  University Press, Cambridge, 1991.

\bibitem{JokP07}
\textsc{S.~Jokar and M.~Pfetsch}, \emph{Exact and approximate sparse solutions
  of underdetermined linear equations}.
\newblock SIAM J. Sci. Comput., To appear, 2007.

\bibitem{LanT85}
\textsc{P.~Lancaster and M.~Tismenetsky}, \emph{The Theory of Matrices},
  Academic Press, New York, 2nd~ed., 1985.

\bibitem{Nat95}
\textsc{B.~K. Natarajan}, \emph{Sparse approximate solutions to linear
  systems}, SIAM J. Comput. \textbf{24}, no.~2 (1995), pp.~227--234.

\bibitem{Str07}
\textsc{G.~Strang}, \emph{Computational Science and Engineering},
  Wellesley-Cambridge Press, 2007.

\bibitem{StrH03}
\textsc{T.~Strohmer and R.~Heath}, \emph{Grassmannian frames with applications
  to coding and communication}, Appl. Comput. Harmon. Anal. \textbf{14}, no.~3
  (2003), pp.~257--275.

\bibitem{Tro04b}
\textsc{J.~A. Tropp}, \emph{Greed is good: algorithmic results for sparse
  approximation}, IEEE Trans. Inform. Theory \textbf{50}, no.~10 (2004),
  pp.~2231--2242.

\end{thebibliography}
\end{document}